\documentclass{article}    
\usepackage{latexsym}
\usepackage{amsbsy}
\usepackage{amssymb}
\usepackage[latin1]{inputenc} 

\usepackage{amsmath}

\title{Diffeomorphism invariant Colombeau algebras. Part II:
Classification}
\author{Michael Grosser\footnote{Electronic mail: michael.grosser@univie.ac.at} \\
  \\{\small Department of Mathematics, University of Vienna}\\
  {\small Strudlhofg. 4, A-1090 Wien, Austria}}
\date{}                                                            
\begin{document}
%
\maketitle
\newcommand{\ep}{\hspace*{\fill}$\Box$}
\newcommand{\eps}{\varepsilon}
\newcommand{\pr}{{\bf Proof. }}
\newcommand{\ms}{\medskip\\}
\newcommand{\cl}{\mbox{\rm cl}}
\newcommand{\g}{\ensuremath{\mathfrak g} }
\newcommand{\gc}{${\cal G}$-complete }
\newcommand{\sa}{\stackrel{\scriptstyle s}{\approx}}
\newcommand{\prol}{\mbox{\rm pr}^{(n)}}
\newcommand{\prolo}{\mbox{\rm pr}^{(1)}}
\newcommand{\deta}{\frac{d}{d \eta}{\Big\vert}_{_{0}}}
\newcommand{\detas}{\frac{d}{d \eta}{\big\vert}_{_{0}}}
\newcommand{\R}{\mathbb R}
\newcommand{\N}{\mathbb N}
\newcommand{\C}{\mathbb C}
\newcommand{\Z}{\mathbb Z}
\newcommand{\K}{\mathbb K}
\newcommand{\sR}{\mathbb R}
\newcommand{\sN}{\mathbb N}
\newcommand{\gK}{{\cal K}_s}   
\newcommand{\gR}{{\cal R}_s}   
\newcommand{\gC}{{\cal C}_s}   
\newcommand{\gKc}{{\cal K}_{sc}}   
\newcommand{\gRc}{{\cal R}_{sc}}   
\newcommand{\Dp}{${\cal D}'$ }                                          
\newcommand{\grn}{${\cal G}(\R^n)$ }
\newcommand{\grp}{${\cal G}(\R^p)$ }
\newcommand{\grq}{${\cal G}(\R^q)$ }
\newcommand{\gt}{${\cal G}_\tau$ }
\newcommand{\gto}{${\cal G}_\tau(\Omega)$ }
\newcommand{\gtrn}{${\cal G}_\tau(\R^n)$ }
\newcommand{\gtrp}{${\cal G}_\tau(\R^p)$ }
\newcommand{\gtrq}{${\cal G}_\tau(\R^q)$ }
\newcommand{\gs}{\ensuremath{{\mathcal G}^s} }
\newcommand{\gso}{\ensuremath{{\mathcal G}^s(\Omega)} }
\newcommand{\gsrn}{\ensuremath{{\mathcal G}^s(\R^n)} }
\newcommand{\gsrp}{\ensuremath{{\mathcal G}^s(\R^p)} }
\newcommand{\gsrq}{\ensuremath{{\mathcal G}^s(\R^q)} }
\newcommand{\gst}{\ensuremath{{\mathcal G}^s_\tau} }
\newcommand{\gsto}{\ensuremath{{\mathcal G}^s_\tau(\Omega)} }
\newcommand{\gstrn}{\ensuremath{{\mathcal G}^s_\tau(\R^n)} }
\newcommand{\gstrp}{\ensuremath{{\mathcal G}^s_\tau(\R^p)} }
\newcommand{\gstrq}{\ensuremath{{\mathcal G}^s_\tau(\R^q)} }
\newcommand{\es}{\ensuremath{{\mathcal E}^s} }
\newcommand{\esm}{\ensuremath{{\mathcal E}^s_M} }
\newcommand{\ns}{\ensuremath{{\mathcal N}^s} }
\newcommand{\est}{\ensuremath{{\mathcal E}^s_\tau} }
\newcommand{\nst}{\ensuremath{{\mathcal N}^s_\tau} }
\newcommand{\ks}{\ensuremath{{\mathcal K}_s} }
\newcommand{\rs}{\ensuremath{{\mathcal R}_s} }
\newcommand{\rsc}{\ensuremath{{\mathcal R}_{sc}} }
\newcommand{\ke}{\ensuremath{{\mathcal K}_e} }
\newcommand{\re}{\ensuremath{{\mathcal R}_e} }
\newcommand{\hsp}{\hspace{-0.3em}}
\newcommand{\gel}{\ensuremath{{\mathcal G}^e} }
\newcommand{\eel}{\ensuremath{{\mathcal E}^e} }
\newcommand{\eelm}{\ensuremath{{\mathcal E}_M^e} }
\newcommand{\nel}{\ensuremath{{\mathcal N}^e} }
\newcommand{\gelt}{\ensuremath{{\mathcal G}^e_\tau} }
\newcommand{\eelt}{\ensuremath{{\mathcal E}^e_\tau} }
\newcommand{\eelmt}{\ensuremath{{\mathcal E}_{M,\tau}^e} }
\newcommand{\nelt}{\ensuremath{{\mathcal N}^e_\tau} }
\newcommand{\interior}[1]{\mathrm{int}(#1)}
\newcommand{\cinfty}{{\cal C}^\infty}
\newcommand{\geltensor}{\ensuremath{{\mathcal G}^{e,\otimes}} }
\newcommand{\catensor}{\ensuremath{{\mathcal A}^\otimes} }
\newcommand{\ketensor}{\ensuremath{{\mathcal K}_{e,\otimes}} }
\newcommand{\gelab}{\ensuremath{{\mathcal G}^{e,\circ}} }
\newcommand{\eelab}{\ensuremath{{\mathcal E}^{e,\circ}} }
\newcommand{\eelmab}{\ensuremath{{\mathcal E}_M^{e,\circ}} }
\newcommand{\nelab}{\ensuremath{{\mathcal N}^{e,\circ}} }
\newcommand{\caab}{\ensuremath{{\mathcal A}^\circ} }
\newcommand{\keab}{\ensuremath{{\mathcal K}_{e,\circ}} }
\newcommand{\uomab}{\ensuremath{U^\circ}(\Om) }
\newcommand{\gtn}{{\cal G}_\tau(\R^n) }
\newtheorem{thr}{\hspace*{-3mm} \bf}[section] 
\newcommand{\bt}{\begin{thr} {\bf Theorem. }}
\newcommand{\et}{\end{thr}}
\newcommand{\bp}{\begin{thr} {\bf Proposition. }}
\newcommand{\bc}{\begin{thr} {\bf Corollary. }}
\newcommand{\blem}{\begin{thr} {\bf Lemma. }}
\newcommand{\bex}{\begin{thr} {\bf Example. }\rm}
\newcommand{\bexs}{\begin{thr} {\bf Examples. }\rm}
\newcommand{\bd}{\begin{thr} {\bf Definition. }}
\newcommand{\beast}{\begin{eqnarray*}}
\newcommand{\eeast}{\end{eqnarray*}}
\newcommand{\wsc}[1]{\overline{#1}^{wsc}}
\newcommand{\todo}[1]{$\clubsuit$\ {\tt #1}\ $\clubsuit$}
\newcommand{\rem}[1]{\vadjust{\rlap{\kern\hsize\thinspace\vbox%
                       to0pt{\hbox{${}_\clubsuit${\small\tt #1}}\vss}}}}
\newcommand{\bappt}{\begin{appthr} {\bf Theorem. }}
\newcommand{\eappt}{\end{appthr}}
\newcommand{\bappp}{\begin{appthr} {\bf Proposition. }}
\newcommand{\bappc}{\begin{appthr} {\bf Corollary. }}
\newcommand{\bappnot}{\begin{appthr} {\bf Notation }}
\newcommand{\caut}[1]{\index{Author}{#1}}
\newcommand{\cind}[1]{\index{Index}{#1}}
\newcommand{\ahat}{\ensuremath{\hat{\mathcal{A}}_0(X)} } 
\newcommand{\atil}{\ensuremath{\tilde{\mathcal{A}}_0(X)} } 
\newcommand{\aqtil}{\ensuremath{\tilde{\mathcal{A}}_q(X)} } 
\newcommand{\amtil}{\ensuremath{\tilde{\mathcal{A}}_m(X)} } 
\newcommand{\ehat}{\ensuremath{\hat{\mathcal{E}}(X)} } 
\newcommand{\emhat}{\ensuremath{\hat{\mathcal{E}}_M(X)} } 
\newcommand{\nhat}{\ensuremath{\hat{\mathcal{N}}(X)} } 
\newcommand{\ghat}{\ensuremath{\hat{\mathcal{G}}(X)} } 
\newcommand{\lhat}{\ensuremath{\hat{L}_\zeta} } 
\newcommand{\comp}{\subset\subset}
\newcommand{\SSS}{{\cal S}}     
\newcommand{\lla}{\langle}\newcommand{\ra}{\rangle} 
\newcommand{\al}{\alpha}
\newcommand{\bet}{\beta} 
\newcommand{\ga}{\gamma}
\newcommand{\Om}{\Omega}\newcommand{\Ga}{\Gamma}\newcommand{\om}{\omega}
\newcommand{\si}{\sigma}\newcommand{\la}{\lambda}
\newcommand{\de}{\delta}
\newcommand{\vphi}{\varphi}\newcommand{\dl}{{\displaystyle \lim_{\eta>0}}\,}
\newcommand{\intl}{\int\limits}\newcommand{\su}{\sum\limits_{i=1}^2}
\newcommand{\D}{{\cal D}}\newcommand{\Vol}{\mathrm{Vol\,}}
\newcommand{\Or}{\mathrm{Or}}\newcommand{\sign}{\mathrm{sign}}
\newcommand{\na}{\nabla}\newcommand{\pa}{\partial}
\newcommand{\ti}{\tilde}\newcommand{\T}{{\cal T}} \newcommand{\G}{{\cal G}}
\newcommand{\DD}{{\cal D}}\newcommand{\X}{{\cal X}}\newcommand{\E}{{\cal E}} 
\newcommand{\CC}{{\cal C}}\newcommand{\vo}{\Vol}
\newcommand{\bat}{\bar t}
\newcommand{\bx}{\bar x}
\newcommand{\by}{\bar y} \newcommand{\bz}{\bar z}\newcommand{\br}{\bar r}
\newcommand{\fr}{\frac{1}}\newcommand{\il}{\int\limits}
\newcommand{\nn}{\nonumber}
\newcommand{\supp}{\mathop{\mathrm{supp}}}

\newcommand{\vp}{\mathrm{vp}\frac{1}{x}}\newcommand{\A}{{\cal A}}
\newcommand{\Ll}{L_{\mathrm{\small loc}}}\newcommand{\Hl}{H_{\mathrm{\small loc}}}
\newcommand{\Lll}{L_{\mathrm{\scriptsize loc}}}
\newcommand{\beq}{ \begin{equation} }\newcommand{\eeq}{\end{equation} }
\newcommand{\bea}{\begin{eqnarray}}\newcommand{\eea}{\end{eqnarray}}
\newcommand{\beas}{\begin{eqnarray*}}\newcommand{\eeas}{\end{eqnarray*}}
\newcommand{\beqs}{\begin{equation*}}\newcommand{\eeqs}{\end{equation*}}
\newcommand{\lb}{\label}\newcommand{\rf}{\ref}
\newcommand{\GL}{\mathrm{GL}}\newcommand{\bfs}{\boldsymbol}
\newcommand{\ben}{\begin{enumerate}}\newcommand{\een}{\end{enumerate}}
\newcommand{\ba}{\begin{array}}\newcommand{\ea}{\end{array}}
\newtheorem{thi}{\hspace*{-1.1mm}}[section]
\newcommand{\bthm}{\begin{thr} {\bf Theorem. }}
\newcommand{\bprop}{\begin{thr} {\bf Proposition. }}
\newcommand{\bcor}{\begin{thr} {\bf Corollary. }}
\newcommand{\bdef}{\begin{thr} {\bf Definition. }}
\newcommand{\brem}{\begin{thr} {\bf Remark. }\rm}
\newcommand{\bth}{\begin{thr}\rm}

\newcommand{\ethi}{\end{thr}}
\newcommand{\tr}{\mathrm{tr}}
\newcommand{\spn}{\mathrm{span}}
\newcommand{\ca}{{\cal A}}
\newcommand{\cb}{{\cal B}}
\newcommand{\cc}{{\cal C}}
\newcommand{\cd}{{\cal D}}
\newcommand{\cee}{{\cal E}}
\newcommand{\cg}{{\cal G}}
\newcommand{\ci}{{\cal I}}
\newcommand{\cn}{{\cal N}}
\newcommand{\css}{{\cal S}}
\newcommand{\ct}{{\cal T}}
\newcommand{\rmd}{\mbox{\rm d}}
\newcommand{\io}{\iota}
\newcommand{\bnot}{\begin{thr} {\bf Notation }}
\newcommand{\lgl}{\langle}
\newcommand{\rgl}{\rangle}
\newcommand{\spp}{\mbox{\rm supp\,}}
\newcommand{\id}{\mathop{\mathrm{id}}}
\newcommand{\pro}{\mathop{\mathrm{pr}}}
\newcommand{\dist}{\mathop{\mathrm{dist}}}
\newcommand{\clb}{\overline{B}_}
\newcommand{\sgn}{\mathop{\mathrm{sgn}}}
\newcommand{\gd}{\ensuremath{{\mathcal G}^d} }
\newcommand{\ed}{\ensuremath{{\mathcal E}^d} }
\newcommand{\edm}{\ensuremath{{\mathcal E}_M^d} }
\newcommand{\nd}{\ensuremath{{\mathcal N}^d} }
\newcommand{\geins}{\ensuremath{{\mathcal G}^1} }
\newcommand{\eeins}{\ensuremath{{\mathcal E}^1} }
\newcommand{\eeinsm}{\ensuremath{{\mathcal E}_M^1} }
\newcommand{\neins}{\ensuremath{{\mathcal N}^1} }
\newcommand{\gzwei}{\ensuremath{{\mathcal G}^2} }
\newcommand{\ezwei}{\ensuremath{{\mathcal E}^2} }
\newcommand{\ezweim}{\ensuremath{{\mathcal E}_M^2} }
\newcommand{\nzwei}{\ensuremath{{\mathcal N}^2} }
\newcommand{\gelnull}{\ensuremath{{\mathcal G}^e_0} }
\newcommand{\eelnull}{\ensuremath{{\mathcal E}^e_0} }
\newcommand{\eelnullm}{\ensuremath{{\mathcal E}_M^e_0} }
\newcommand{\nelnull}{\ensuremath{{\mathcal N}^e_0} }
\newcommand{\supl}{\sup\limits}

\begin{abstract}     
This contribution presents a comprehensive analysis of
Colombeau \!\!\!\hbox{(-type)} algebras in the range between the 
diffeomorphism invariant algebra $\gd =
\edm\big/\nd$ introduced in Part I (see \cite{found})
and Colombeau's original algebra $\mathcal{G}^e$ introduced in
\cite{c2}. Along the way, it provides several classification
results (again see
\cite{found})
which are indispensable for obtaining an intrinsic
description of a (full) Colombeau algebra on a manifold (\cite{vim}). The
latter will be the focus of Part III of this series of
contributions.
\end{abstract}

\section{Introduction}\lb{introII}\lb{calc2}
This contribution continues the first in a series of three
(Parts I and III also in this volume)
by analyzing diffeomorphism invariant Colombeau
algebras from a broader point of view.
We will use freely notations and results from Part I;
for details see \cite{found}.

The main result of Section \rf{cond0}
below allows for considerably simplifying the definition
of the null ideal: Indeed, it dispenses with taking into
account the derivatives of the representative being tested.
This applies to virtually all versions of Colombeau algebras.
In Section \rf{ex} we 
show that
the diffeomorphism invariant algebra $\cg^d(\Om)$
of \cite{Jel} resp.\ \cite{found}
(see Section 3 of Part I)
is not injectively included in
the Colombeau algebra $\cg^e(\Om)$ of \cite{c2}
by constructing two counterexamples.
Section \rf{spec} develops a framework allowing to classify
the range of algebras 
which can be positioned between $\cg^d(\Om)$ and
(the smooth version of) $\cg^e(\Om)$.
In particular, we are going to
determine the minimal extent to which
the definition of the  algebra
introduced by J.~F.~Colombeau and A.~Meril in
\cite{CM}
has to be modified to obtain diffeomorphism invariance.
This leads to the construction of the (diffeomorphism
invariant) Colombeau algebra  $\cg^2(\Om)$ which is closer to the
algebra  of \cite{CM}
than the algebra $\cg^d(\Om)$.
Certain classification
results of Section
\rf{spec} are essential for obtaining an intrinsic
description of Colombeau algebras on manifolds (see Part III resp.\
\cite{vim}).

Both the counterexamples to be constructed in Section \rf{ex}
will take the
form of infinite series, being absolutely convergent in each
derivative. Thus we need a theorem
guaranteeing the completeness of
$\cee(\Om)=\cc^\infty(U(\Om),\C)\equiv
    \cc^\infty(\ca_0(\Om)\times\Om,\C)$
with respect to the corresponding topology.
To this end, let $E,F$ be locally convex spaces
and $U$ an open subset of $E$. If $f:U\to F$ is
smooth, its
$n$-th differential $\rmd^n\!f$ belongs to
$\cc^\infty(U,L^n(E^n,F))$ where $L^n(E^n,F)$ denotes the space
$L(E,\dots,E;F)$ of $n$-linear bounded maps from
$E\times\dots\times E$ ($n$ factors) into $F$.
(For $n=0$, set $L^n(E^n,F):=F$.)
On $\cc^\infty(U,L^n(E^n,F))$, let $\tau_{cb}^n$ denote
the topology of uniform $F$-convergence on subsets of the form $K\times B$
where $K$ is a compact subset of $U$ and $B$ is bounded in
$E^n=E\times\dots\times E$. Let $\cc^\infty(U,F)$ carry the initial
(locally convex) topology $\tau^\infty$ induced by the family
$(\rmd^n,\cc^\infty(U,L^n(E^n,F)),\tau_{cb}^n)_{n\ge0}$, i.e., the
topology of uniform convergence of all derivatives
(that is to say, differentials) on sets $K\times B$ as above.
Note that on $\cc^\infty(\R,F)$, $\tau^\infty$ 
is just the usual Fr\'echet topology of compact convergence in all
derivatives. For the proof of the following theorem, see
\cite{found}.

\bt\lb{cinftcomp}
Let $E,F$ be locally convex spaces, assume $F$ to be complete
and let $U$ be an open subset of $E$. Then $\cc^\infty(U,F)$
is complete with respect to the topology
$\tau^\infty$ of uniform $F$-convergence of all differentials
on subsets of the form $K\times B$
where $K$ is a compact subset of $U$ and $B$ is bounded in
the appropriate product $E^n=E\times\dots\times E$. Moreover,
for each $p\in\N$, the operator
$\rmd^p:\cc^\infty(U,F)\to\cc^\infty(U,L^p(E^p,F))$ is continuous
if both the domain and the range space carry
the respective topology $\tau^\infty$.
\et

Now let $U$ denote a (non-empty) open subset
of a closed affine subspace $E_1$ of some locally convex space $E$,
$E_0$ the linear subspace parallel to $E_1$ and $F$ a complete
locally convex space. {\it Mutatis mutandis},
\rf{cinftcomp} is valid also in this slightly more general
situation. The vectors
$v_1,\dots,v_n$ to be plugged into $\rmd^n\!f(x)$ now have to be
taken from $E_0$, as well as $B$ has to denote a bounded subset
of $E_0^n$.

In the following, we will abbreviate $R\circ S^{(\eps)}$ as
$R_\eps$, throughout. Terms of the form $\pa^\al\rmd_1^kR_\eps$
always are to be read as $\pa^\al\rmd_1^k(R_\eps)$.

\section{A simple condition equivalent to negligibility}\lb{cond0}
The principal part of this section refers to $\gd$. However, in the
concluding remarks we will indicate that the main result is true for
virtually all types of Colombeau algebras.

Th.\ 18\,(2$^\circ$) of \cite{Jel}
gives a condition equivalent to negligibility
replacing the term $\pa^\al(R(S_\eps\phi(\eps,x),x))$ occurring in
the definition (3.1 in Part I resp.\ Def.\ 7.3 in \cite{found}) by 
$(\pa^\al\rmd_1^k R_\eps)(\vphi,x)(\psi_1,\dots,\psi_k)$.
(The analogue of this theorem
for the case of moderateness can be looked up as 3.2
in Part I.)
Moreover, Th.\ 18\,(1$^\circ$) of \cite{Jel} shows that
we still get a condition 
equivalent to $R\in\cn(\Om)$ if we
simply omit the differential with respect to the first variable
$\vphi$ from
(2$^\circ$), provided $R$ is assumed to be
moderate. In the following, we are going to show that
a further simplification is possible which might seem rather
drastic at first glance: It is not even necessary to
consider partial derivatives with respect to $x\in\Om$.
In order to facilitate comparing the conditions mentioned 
so far we include all of them in the following theorem,
though only $(0^\circ)$ is new.
\bt\lb{thcond0}
For $R\in\cee_M(\Om)$,
each of the following conditions is equivalent
to $R\in\cn(\Om)$:\medskip\\
$(0^\circ)$\quad
$\forall K\subset\subset\Om\  
   \forall n\in\N\ \exists q\in\N
   \ \forall B\,(\mbox{bounded})\,\subseteq \cd(\R^s)$:
$$R_\eps(\vphi,x)=O(\eps^n)
   \qquad\qquad (\eps\to0).$$
$(1^\circ)$\quad
$\forall K\subset\subset\Om\ \forall\al\in\N_0^s\ 
   \forall n\in\N\ \exists q\in\N
   \ \forall B\,(\mbox{bounded})\,\subseteq \cd(\R^s)$:
$$\pa^\al R_\eps(\vphi,x)=O(\eps^n)
   \qquad\qquad (\eps\to0).$$
$(2^\circ)$\quad
$\forall K\subset\subset\Om\ \forall\al\in\N_0^s\ 
   \forall k\in\N_0\ \forall n\in\N\ \exists q\in\N
   \ \forall B\,(\mbox{bounded})\,\subseteq \cd(\R^s)$:
$$\pa^\al\rmd_1^k R_\eps(\vphi,x)(\psi_1,\dots,\psi_k)=O(\eps^n)
   \qquad\qquad (\eps\to0).$$
In each of the preceding conditions, the estimate is to be
understood as to hold uniformly with respect to
$x\in K$, $\vphi\in B\cap\ca_q(\R^s)$ $((1^\circ)$,$(2^\circ))$,
   $\psi_1,\dots,\psi_k\in B\cap\ca_{q0}(\R^s)$ $((2^\circ))$.
\et
\pr
To highlight the part of the theorem which is new as compared to
Th.\ 18 of \cite{Jel} we present the proof of
$(0^\circ)\Rightarrow(1^\circ)$.
To this end, we will show, assuming $R\in\cee_M(\Om)$ to satisfy
$(0^\circ)$, that $R$ satisfies $(1^\circ)$ for
$\al:=e_i$, i.e., $\pa^\al=\pa_i$ ($i=1,\dots,s$)
and that, in addition, $\pa_i R$ again is moderate and satisfies 
$(0^\circ)$. Then it will follow by induction that $(1^\circ)$
holds for all $\al\in\N_0^s$.

So suppose $R\in\cee_M(\Om)$ to satisfy $(0^\circ)$
and let $K\subset\subset\Om$ and
$n\in\N$ be given. For $\de:=\min(1,\dist(K,\pa\Om))$,
set $L:=K+\clb{\frac{\de}{2}}(0)$.
Then $K\subset\subset L\subset\subset\Om$.
Now by moderateness of $R$ and Th.\ 17 of \cite{Jel}
(3.2 of Part I), choose $N\in\N$
such that for every bounded subset $B$ of $\cd(\R^s)$ the relation
$\pa_i^2R_\eps(\vphi,x)=O(\eps^{-N})$ as $\eps\to0$
holds,
uniformly for $x\in L$, $\vphi\in B\cap\ca_0(\R^s)$. Next,
by the assumption of $(0^\circ)$ to hold for $R$, choose
$q\in\N$ such that, again for every bounded subset $B$ of
$\cd(\R^s)$, we have
$R_\eps(\vphi,x)=O(\eps^{2n+N})$ as $\eps\to0$,
uniformly for $x\in L$, $\vphi\in B\cap\ca_q(\R^s)$.
Now suppose a bounded subset $B$ of $\cd(\R^s)$ to be given;
let $\vphi\in B\cap\ca_q(\R^s)$, $x\in K$ and
$0<\eps<\frac{\de}{2}$; hence $x+\eps^{n+N}e_i\in L$.
By Taylor's Theorem, we conclude
(to be precise, separately for the real and imaginary part of $R$)
\beas
R_\eps(\vphi,x+\eps^{n+N}e_i)
   &=&R_\eps(\vphi,x)+\pa_iR_\eps(\vphi,x)\eps^{n+N}+
      \frac{1}{2}\pa_i^2R_\eps(\vphi,x_\theta)\eps^{2n+2N}
\eeas
where $x_\theta=x+\theta\eps^{n+N}e_i$ for some $\theta\in(0,1)$;
note that also $x_\theta\in L$.
Consequently,
\beas
\pa_iR_\eps(\vphi,x)
   &=&
       \underbrace
         {\left(R_\eps(\vphi,x+\eps^{n+N}e_i)
          -R_\eps(\vphi,x)\right)}
         _{O(\eps^{2n+N})} \eps^{-n-N}
       -\underbrace
         {\frac{1}{2}\pa_i^2R_\eps(\vphi,x_\theta)}
         _{O(\eps^{-N})}
       \eps^{n+N},
       \eeas
uniformly for $\vphi\in B\cap\ca_q(\R^s)$, $x\in K$.
Having demonstrated $\pa_iR_\eps(\vphi,x)=O(\eps^n)$
for all $i=1,\dots,s$, observe that
$\pa_i(R_\eps)=(\pa_iR)_\eps$. Therefore,
$\pa_iR$ again satisfies $(0^\circ)$.
According to Th.\ 7.10 of \cite{found} (which is non-trivial,
see the discussion in Section 7.3 of \cite{found}),
$\pa_iR$ is also moderate .
By the remark made above, this
completes the proof.
\ep

The reader acquainted with E.~Landau's paper \cite{Lan} will
easily recognize the method employed therein to form the basis of
the preceding proof.

The seemingly technical difference
between $(0^\circ)$ and the remaining conditions
has decisive effects on
applications: For example, if the uniqueness of a
solution of a differential equation is to be shown
one supposes
$R_1,R_2$ to be representatives of solutions. Note that this
includes the assumption that $R_1,R_2\in\cee_M(\Om)$, hence
\rf{thcond0} may be applied. For $[R_1]=[R_2]$ in $\cg(\Om)$
we have to show that $R:=R_1-R_2\in\cn(\Om)$. Now it suffices to check
condition $(0^\circ)$ rather than $(1^\circ)$ (resp.\ $(2^\circ)$
resp.\ the original definition of $R\in\cn(\Om)$),
i.e., 
there is no need
to analyze the behaviour of any derivative of $R$.

The part of \ref{thcond0}
saying that for moderate functions (the appropriate analog of)
condition $(0^\circ)$ is equivalent to 
negligibility
applies to virtually all versions of
Colombeau algebras of practical importance,
in particular, to the following:
\begin{itemize}
\item
For the special algebra as defined,
e.g., in \cite{MOBook}, p.~109, just replace the term
$R_\eps(\vphi,x)$
in condition $(0^\circ)$ by $u_\eps(x)$.
\item
For the classical full Colombeau algebra
of \cite{c2} simply drop the uniformity requirement
concerning $\vphi$ from $(0^\circ)$.
\item
For the diffeomorphism invariant Colombeau algebra
$\cg^2(\Om)$ to be introduced in Section \rf{g2},
the corresponding result is stated as Th.\ 17.9
in \cite{found}.
\item
For the special algebra on smooth manifolds
the corresponding result follows from the local
characterization of generalized functions
(see \cite{Stdiss}, 4.4).
\item
The latter also applies to the intrinsically defined
full Colombeau algebra on manifolds (\cite{vim}, Cor.\ 4.5).
\end{itemize}
In the first and second of these four instances, the
respective proofs are obtained by appropriately slimming down
the proof of
\rf{thcond0}.

\section{Non-injectivity of the canonical homomorphism  from
$\cg^d(\Om)$ into $\cg^e(\Om)$}\lb{ex}

For every open subset $\Om$ of $\R^s$, there is a canonical
algebra homomorphism $\Phi$ from the diffeomorphism invariant
Colombeau algebra $\cg^d(\Om)$ of \cite{Jel} (see Section
3 of Part I) to the ``classical'' (full) Colombeau algebra
$\cg^e(\Om)$ introduced in \cite{c2}, 1.2.2 (see Section 1 of Part
I). By constructing suitable (coun\-ter)ex\-am\-ples,
we are going to show that $\Phi$ is not injective in general.

By superscripts $d,e$ we will distinguish between ingredients
for constructing 
$\cg^d$ resp.\ $\cg^e$.
As in Section 3 of Part I we will use the C-formalism
also in the present context.
To see that
$\cee^d$ is a subset of $\cee^e$
we have to pass from C-representatives to
J-representatives: Smoothness of $R^d\in\cee^d$, by definition,
is equivalent to smoothness of
$(T^*)^{-1}R^d\in\cc^\infty(\ca_0(\Om)\times\Om)$ while
for $R^e\in\cee^e$, smoothness of $x\mapsto R^e(\vphi,x)$
is equivalent to smoothness of
$x\mapsto(T^*)^{-1}R^e(\vphi(.-x),x)$. From this it is clear that
$\cee^d\subseteq\cee^e$.
Moreover, we obtain $\cee^d_M\subseteq\cee^e_M$ and
$\cn^d\subseteq\cn^e$.
This follows easily by
inspecting the corresponding definitions.

Thus we obtain a canonical map $\Phi:\cg^d(\Om)\to\cg^e(\Om)$
which is an algebra homomorphism
respecting the embeddings of $\cd'(\Om)$ and differentiation.

\brem (i)
Colombeau's original
construction in 1.2.2 of \cite{c2} produces a full algebra
$\gel_1(\Om)$ differing slightly from $\gel(\Om)$ used above.
$\gel_1(\Om)$ \linebreak is obtained on the basis of
$U_1(\Om):=T^{-1}(\ca_1(\Om)\times\Om)$ rather than
$U(\Om)=T^{-1}(\ca_0(\Om)\times\Om)$.
The
restriction operator $\Phi_0$ maps $\cee^d$ into
$\cee^e_1$, $\cee^d_M$ into $\cee^e_{1,M}$ and $\cn^d$ into
$\cn^e_1$,
respectively. The canonical map $\Phi_1:\gd\to\gel_1$ induced by $\Phi_0$
acts on representatives as
restriction from $T^{-1}(\ca_0(\Om)\times\Om)$ to
$T^{-1}(\ca_1(\Om)\times\Om)$.\smallskip\\
(ii) The counterexamples to be constructed below will settle the
question of injectivity not only of
$\Phi:\cg^d\to\gel$ but also of
$\Phi_1:\cg^d\to\gel_1$:
$\Phi_1$ is injective if and only if $\Phi$ is,
due to the canonical map
$\Psi:\gel\to\gel_1$ being injective.
\et

In the following, we will define maps $P,Q:U(\R)\to\C$ each of
which satisfies the following conditions (i)--(iv),
thereby providing a counterexample
to the conjecture of the canonical map $\Phi$ being injective.
\begin{itemize}
\item[(i)] $R\in\cee^d$, i.e., $R$ has to be smooth;
\item[(ii)] $R\in\cee_M^d$,
\item[(iii)]$R\notin\cn^d$,
\item[(iv)] $R\in\cn^e$.
\end{itemize}
Let $s:=1$, $\Om:=\R$. As a prerequisite
we introduce the following
notation:
$$\ba{lrcll}
&\lgl\vphi|\vphi\rgl&:=&\int\vphi(\xi)\overline{\vphi(\xi)}\,d\xi
      \qquad\qquad &(\vphi\in\cd(\R))\medskip\\
v_k\in\cd'(\R):\qquad&\lgl v_k,\vphi\rgl&:=&\int\xi^k\vphi(\xi)\,d\xi
      \qquad\qquad&(\vphi\in\cd(\R),\ k\in\N_0)\medskip\\
v_\frac{1}{2}\kern-2pt\in\cd'(\R):\qquad&\lgl v_\frac{1}{2},
                     \vphi\rgl&:=&\int|\xi|^\frac{1}{2}\vphi(\xi)\,d\xi
      \qquad\qquad&(\vphi\in\cd(\R))\medskip\\
&v(\vphi)&:=&\lgl\vphi|\vphi\rgl^\frac{1}{2}
                              \lgl v_\frac{1}{2},\vphi\rgl
      \qquad\qquad&(\vphi\in\cd(\R))\medskip\\
&g(x)&:=&\frac{x}{1+x^2}&(x\in\R)\medskip\\
&e(x)&:=&\begin{cases} 
             \exp(-\frac{1}{x})\quad&(x>0)\\
             0                     &(x\le0)
          \end{cases} &(x\in\R)\medskip\\
&\ga_k&:=&k+\frac{1}{k}&(k\in\N).
\ea$$
Finally, choose an (even) function
$\si\in\cd(\R)$ satisfying $0\le\si\le1$,
$\si(x)\equiv1$ for $|x|\le\frac{1}{2}$,
$\si(x)\equiv0$ for $|x|\ge\frac{3}{2}$ and set
$$h_k(x):=\si(x)\cdot 2g(x)+(1-\si(x))\cdot
\sgn(x)\cdot|2g(x)|^{\ga_k}\qquad\qquad(x\in\R,
 \ k\in\N).$$

Apart from abbreviating $R\circ S^{(\eps)}=R\circ(S_\eps\times\id)$
as $R_\eps$ for any function $R$ defined on $\ca_0(\R)\times\R$,
we also will write $R_\eps$ for $R\circ S_\eps$ if $R$ is defined on
$\ca_0(\R)$.
\bd \lb{defpq}
Let $\vphi\in\ca_0(\R)$, $x\in\R$ and set
\beas
P(\vphi,x)&:=&\sum_{k=1}^{\infty}\frac{1}{k!}\cdot
   g\big(\lgl\vphi|\vphi\rgl^{\ga_k} e(v(\vphi))\big)\cdot
   \lgl\vphi|\vphi\rgl^{\ga_k}\cdot
   \lgl v_k,\vphi\rgl,\\
Q(\vphi,x)&:=&\sum_{k=1}^{\infty}\frac{1}{k!}\cdot
   h_k\big(\lgl\vphi|\vphi\rgl^{\frac{3}{2}}\, \lgl v_{\frac{1}{2}},
                                   \vphi\rgl\big)\cdot
   \lgl\vphi|\vphi\rgl^{\ga_k}\cdot
   \lgl v_k,\vphi\rgl.
\eeas
\et
Hence $P$ and $Q$, in fact, only depend on $\vphi$. 
Explicitly, $P$ is given by
\beas \lefteqn{
P(\vphi,x)=}\\
    &&\sum_{k=1}^{\infty}\frac{1}{k!}\cdot
   \frac{\left(\int\vphi(\xi)\overline{\vphi(\xi)}\,d\xi\right)^{{k+\frac{1}{k}}}
   \exp\left(-\frac{1}{\left(\int\vphi(\xi)\overline{\vphi(\xi)}\,d\xi\right)^\frac{1}{2}
   \left(\int|\xi|^\frac{1}{2}\vphi(\xi)\,d\xi\right)}\right)}
   {1+\left(\left(\int\vphi(\xi)\overline{\vphi(\xi)}\,d\xi\right)^{{k+\frac{1}{k}}}
   \exp\left(-\frac{1}{\left(\int\vphi(\xi)\overline{\vphi(\xi)}\,d\xi\right)^\frac{1}{2}
   \left(\int|\xi|^\frac{1}{2}\vphi(\xi)\,d\xi\right)}\right)\right
   )^2}.
\eeas

It can be shown 
that the series for both $P$ and $Q$ converge uniformly on
bounded subsets of $\ca_0(\R)$, rendering $P$ and $Q$ well-defined
by \ref{cinftcomp}. For the proof of claims (i)--(iv) above we
refer to \cite{found}.

The reader might ask if it is indeed necessary to come up with
counterexamples as complicated as $P$ and $Q$ certainly are.
The author doubts
that easier ones might be possible. This view is based on reflecting
on the r\^oles each of the three factors constituting
a single term of the series
(for $P$, say) in fact has to play:
\begin{itemize}
\item $\lgl v_k,\vphi\rgl$ distinguishes between the spaces
      $\ca_q(\R)$; this is crucial for the negligibility
      properties.
\item $\lgl\vphi|\vphi\rgl^{\ga_k}=\lgl\vphi|\vphi\rgl^k\cdot
      \lgl\vphi|\vphi\rgl^{\frac{1}{k}}$, on the one hand,
      after scaling of $\vphi$
      compensates for the factor $\eps^k$ generated by scaling
      $\vphi$ in $\lgl v_k,\vphi\rgl$. On the other hand, it
      introduces a factor $\eps^{-\frac{1}{k}}$ making the
      first non-vanishing term of the series the dominant one
      as $\eps\to0$.
\item $g(\lgl\vphi|\vphi\rgl^{\ga_k} e(\lgl v,\vphi\rgl))$
      allows the pointwise vs.\ uniformly distinction being
      necessary to obtain $P\notin\cn^d$, $P\in\cn^e$.
      Though
      $g(\lgl\vphi|\vphi\rgl^{\ga_k}\lgl v,\vphi\rgl)$ would suffice
      to achieve the latter, this alternative choice
      for the argument of $g$ would produce, via the chain rule,
      a factor $\eps^{-n(k+\frac{1}{k})}$ in the $k$-th term of
      $\rmd^n\!P_\eps$ which would be disastrous for the moderateness of
      $P$. The function $e$ (together with $\eps^{-\ga_k}$ in the
      argument of $g$) suppressing this
      unwanted factor, $P$ becomes moderate in the end.
\end{itemize}
Similar arguments apply to $Q$.
\section{Classification of smooth Colombeau algebras between $\cg^d(\Om)$
and $\cg^e(\Om)$}\lb{spec}      \lb{g2}
%
Apart from $\cg^e(\Om)$, all algebras to be considered in this
section
have $\cc^\infty(U(\Om))$ resp.\
$\cc^\infty(\ca_0(\Om)\times\Om)$ as their basic space. In
particular, they are smooth algebras in the sense that
representatives $R$ have to be smooth also with respect to $\vphi$.
The term ``test object'' will always refer to some element of
$\cc^\infty_b(I \times\Om,\ca_0(\R^s))$.
\bd
Let $q\in\N$. A
function $\phi:I \to\cd(\R^s)$ (possibly 
depending also on other arguments, e.g., on $x\in\Om$) 
is said to have {\rm vanishing moments of order $q$}
if $\int\xi^\al\phi(\eps)(\xi)\,d\xi=0$ for all $\al\in\N_0^s$
with $1\le|\al|\le q$. It is said to have
{\rm asymptotically vanishing moments of order $q$}
if $\int\xi^\al\phi(\eps)(\xi)\,d\xi=O(\eps^q)$
for all $\al\in\N_0^s$ with $1\le|\al|\le q$. To which extent
this estimate is assumed to hold uniformly with respect
to, e.g., $x\in\Om$ has to be specified separately
(see below).
\et

To obtain a classification of
Colombeau algebras lying in the range
between $\cg^d(\Om)$ and (the smooth version of)
$\cg^e(\Om)$ we introduce symbols  of the forms
$[\mathrm{p}]$, $[\mathrm{M}]$, $[\mathrm{p,M}]$
where p refers to the parameters and M to the moment
properties of a test object. p, being one of c,
$\eps$, $\eps x$ denotes test objects of the form
$\vphi$ (``constant''), $\phi(\eps)$ and $\phi(\eps,x)$,
respectively. M, on the other hand, can take the values 0,A,V,
corresponding to $\ca_0(\R^n)$, asymptotically vanishing moments
and $\ca_q(\R^n)$, respectively. $[\mathrm{A}]$ only applies to
parametrization type $[\eps]$. For test objects of type
$[\eps x]$, we distinguish the following uniformity requirements
concerning
asymptotically vanishing moments:

$[\mathrm{A}_\mathrm{l}]$:\hspace{.5em} uniformly on the particular $K\subset\subset\Om$
(``locally'');

$[\mathrm{A}_\mathrm{g}]$:\hspace{.4em} uniformly on each $L\subset\subset\Om$
        (``globally'');

$[\mathrm{A}_\mathrm{l}^\infty]$: all derivatives  $\pa_x^\al\phi(\eps,x)$
        uniformly on the particular $K\subset\subset\Om$;

$[\mathrm{A}_\mathrm{g}^\infty]$: all derivatives  $\pa_x^\al\phi(\eps,x)$
        uniformly on each $L\subset\subset\Om$.

\noindent
Here, ``on the particular $K\comp\Om$'' is to be read as
``on the particular $K\comp\Om$ on which $R$ is being tested''.
If this compact set $K$ and/or the
order $q$ of the (asymptotic) vanishing of moments is to be
specified, $K$ resp.\ $q$ will be put as subscript(s)
to the corresponding
A-symbol, e.g., $[\mathrm{A}_\mathrm{l}]_{K,q}$.
If in $[\mathrm{p},\mathrm{M}]$ $\mathrm{M}$ is one of the A-symbols
then $\mathrm{p}=\eps$ resp.\ $\mathrm{p}=\eps x$,
being redundant, will be omitted
frequently.

If $[X]$ and $[Y]$ are chosen from the set of the eleven types
such that
$\cee_M[X]\subseteq\cee_M[Y]$) and if,
in addition, $[Y]$ is one of the types $[\mathrm{A}]$ or
$[\mathrm{V}]$ then it easily checked that $\cee_M[X]$ is an algebra
containing $\cn[Y]\cap\cee_M[X]$ as an ideal.
Consequently, $\cee_M[X]\big/(\cn[Y]\cap\cee_M[X])$ is an
algebra. We shall refer to algebras arising in this way by the term
``Colombeau-type algebras''. Altogether there are 46 admissible
choices of pairs $[X],[Y]$. In the following definition, we will
specify eleven algebras of this kind, one for each type of
moderateness. These will be the only ones we are to deal with in
the sequel.
Each of the
remaining Colombeau-type algebras can be obtained as some subalgebra
or some quotient algebra of one of them. Note, however, that the
collection of these eleven algebras is not minimal in this respect
(see Th.\ 17.10 in \cite{found}).
\bd\lb{colalg}
If $[X]$ is one of the types $[\mathrm{V}]$ or $[\mathrm{A}]$
define
$$\cg[X]:=\cee_M[X]\big/\cn[X];$$
for types $[0]$ define
\beas
  \cg[\eps x,0]&:=&\cee_M[\eps x,0]\big/
    \big(\cn[\eps x,\mathrm{A}_\mathrm{l}^\infty]\cap\cee_M[\eps
    x,0]\big),\\
  \cg[\eps,0]&:=&\cee_M[\eps,0]\big/\big(\cn[\eps,\mathrm{A}]
    \cap\cee_M[\eps,0]\big),\\
  \cg[\mathrm{c},0]&:=&\cee_M[\mathrm{c},0]\big/\big(\cn[\mathrm{c},
  \mathrm{V}]\cap\cee_M[\mathrm{c},0]\big).
\eeas
\et
%
We will refer to $\cg[X]$ also by ``the algebra of type $[X]$''.
The open set $\Om$ is omitted from the notation.
Denoting by $\cg^e_0(\Om)$ the ``smooth part'' of $\cg^e(\Om)$,
i.e., the subalgebra formed by all members having a smooth representative
$R\in\cc^\infty(U^e(\Om))$, it is easy to see that
$\cg^e_0(\Om)=\cg[\mathrm{c},\mathrm{V}]$.
$\cg^1(\Om)$ obviously is equal to $\cg[\eps,\mathrm{A}]$; the
algebra $\cg^2(\Om)$ to be
discussed below is obtained as $\cg[\eps x,\mathrm{A}_\mathrm{g}^\infty]$.
$\cg^d(\Om)$, finally, is given as $\cg[\eps x,0]$.
Observe that according to
Th.\ 7.9 of \cite{found} (3.4 in Part I), $\cn[\eps
x,\mathrm{A}_\mathrm{l}^\infty]$ can be
replaced by $\cn[\eps x,\mathrm{V}]$ in the definition of
$\cg[\eps x,0]$.

Cor.\ 16.8 of \cite{found} shows that
test objects of types $[\mathrm{A}_{\mathrm{g}}]$ and
$[\mathrm{A}^\infty_{\mathrm{g}}]$, respectively,
give rise to the same moderate resp.\ negligible functions.
Moreover, by Cor.\ 17.6 of \cite{found} also
test objects of type $[\mathrm{A}_\mathrm{l}^\infty]$ lead to the same
respective notions of moderateness and negligibility
as test objects of type
$[\mathrm{A}^\infty_{\mathrm{g}}]$ do.
This actually leaves us with nine possibly different
algebras. 

As to the diagram formed by
the canonical homomorphisms between these nine algebras,
note that there is no
such mapping from $\cg^d(\Om)=\cg[\eps x,0]$ into
$\cg[\eps x,\mathrm{A}_\mathrm{l}]$
since
$\cn[\eps x,\mathrm{A}_\mathrm{l}]\cap\cee_M[\eps x,0]$---not
containing any of the functions \linebreak
$R(\vphi,x):=\int\xi^\bet\vphi(\xi)\,d\xi$---is
strictly smaller than
$\cn[\eps x,\mathrm{A}_\mathrm{l}^\infty]\cap\cee_M[\eps x,0]$.
We do have canonical homomorphisms, however, both from
$\cg^d(\Om)=\cg[\eps x,0]$ and from $\cg[\eps x,\mathrm{A}_\mathrm{l}]$
into $\cg^2(\Om)=\cg[\eps x,\mathrm{A}_\mathrm{g}^\infty]$.
So we finally arrive at
$$\begin{array}{ccccccc}
{}&  &\cg[\eps x,0]&\to&\cg[\eps,0]&\to&\cg[\mathrm{c},0]\\
&&\downarrow&&\downarrow&&\\
{}\cg[\eps x,\mathrm{A}_\mathrm{l}]&\to&\cg[\eps x,\mathrm{A}^\infty_{\mathrm{g}}]
   &\to&\cg[\eps,\mathrm{A}]&&\downarrow\\
&&\downarrow&&\downarrow&&\\
&&\cg[\eps x,\mathrm{V}]&\to&\cg[\eps,\mathrm{V}]&\to&\cg[\mathrm{c,V}]
\end{array}$$

By the methods employed in \cite{found} one can show that each of
the nine algebras occurring in the diagram (injectively) contains
$\cd'(\Om)$ via $\io$; with one exception (namely, 
$\cg[\eps x,\mathrm{A}_\mathrm{l}]$; cf.\ Ex.\ 7.7 in \cite{found})
the restriction of $\io:\cd'\to\cg$ to $\cinfty$ coincides with
$\si:\cinfty\to\cg$, implying that $\io$ preserves the product of
smooth functions, see \cite{found} for details and proofs.
Moreover, for each type $[X]$ except
$[\eps x,\mathrm{A}_\mathrm{l}]$, $\cee_M$ and $\cn$ are invariant
under differentiation, thus rendering $\cg[X]$ a
differential algebra. Concerning diffeomorphism invariance,
finally, one can
show that $\gd=\cg[\eps x,0]$,
$\cg^2:=\cg[\eps x,\mathrm{A}^\infty_{\mathrm{g}}]$
and $\cg[\eps x,\mathrm{A}_\mathrm{l}]$ in fact share this
property,
yet neither of the remaining six algebras does. For the proofs of
these statements we refer to Chapter 17 of \cite{found}.
$\cg^2:=\cg[\eps x,\mathrm{A}^\infty_{\mathrm{g}}]$
turns out to be the most delicate case in the technical respect.

Summarizing, we obtain that
$\cg^d(\Om)$ and $\cg^2(\Om)$ are the only diffeomorphism
invariant Colombeau algebras among the eleven (resp.\ nine)
algebras defined in \rf{colalg}.

The algebra $\cg^2(\Om)$ 
of type $[\eps x,\mathrm{A}_\mathrm{g}^\infty]$
can be viewed as resulting
from the algebra $\cg^1(\Om)=\cg[\eps,\mathrm{A}]$
of \cite{CM} by applying the minimal modification necessary to
obtain diffeomorphism invariance.

The fact that all three types $[\mathrm{A}_\mathrm{g}^\infty]$,
$[\mathrm{A}_\mathrm{g}]$ and $[\mathrm{A}_\mathrm{l}^\infty]$
give rise to the same notions of moderateness resp.\ negligibility,
hence to the same Colombeau algebra, constitutes one of
the key ingredients for obtaining an intrinsic description of
the algebra $\cg^d$ on manifolds: 
The property of a test object
living on the manifold to have asymptotically vanishing moments can
be formulated in intrinsic terms, indeed
(see \cite{vim}, Def.\ 3.5 resp.\ Part III); yet it would be virtually
unmanageable to deal with the latter property also for derivatives
of this test object, which, of course, are to be understood in this
general case as
appropriate Lie derivatives with respect to smooth vector fields.
Now Cors.\ 16.8 and 17.6 of \cite{found}
allow to dispense with derivatives of test objects as regards the
asymptotic vanishing of the moments, provided all
$K\subset\subset\Om$ are taken into account (\cite{vim}, Cor.\ 4.5).


\begin{thebibliography}{99}
%
%
\bibitem[{Col}85]{c2}
{Colombeau, J.~F.}
\newblock {\em Elementary Introduction to New Generalized Functions}.
\newblock North Holland, Amsterdam, 1985.

\bibitem[{Col}94]{CM}
{Colombeau, J.~F., Meril, A.}
\newblock Generalized functions and multiplication of distributions on
  {${\mathcal C}^\infty$} manifolds.
\newblock {\em J.~Math.~Anal.~Appl.}, {\bf 186}:357--364, 1994.
%
\bibitem[{Gro}99]{vim}
{Grosser, M., Kunzinger, M., Steinbauer, R., Vickers, J.}
\newblock A global theory of algebras of generalized functions.
\newblock {\em Preprint (available electronically at {\tt
  http://arXiv.org/abs/math.FA/9912216})}, 1999.

\bibitem[{Gro}01]{found}
{Grosser, M., Farkas, E., Kunzinger, M., Steinbauer, R.}
On the foundations of nonlinear generalized functions {I}, {II}.
{\em Mem.~Am.~Math.~Soc., to appear (available electronically at {\tt
  http:}\,{\tt //arXiv.org/abs/math.FA/9912214, 9912215})}, 2001.

\bibitem[{Jel}99]{Jel}
{Jel\'\i nek, J.}
\newblock An intrinsic definition of the {C}olombeau generalized functions.
\newblock {\em Comment.~Math.~Univ.~Carolinae}, {\bf 40}:71--95, 1999.

%
%
%
%
\bibitem[{Lan}14]{Lan}
{Landau, E.}
\newblock Einige {U}ngleichungen f\"ur zweimal differentiierbare {F}unktionen.
\newblock {\em Proc.~London Math.~Soc. Ser.~2}, {\bf 13}:43--49, 1913--1914.
%
\bibitem[{Obe}92]{MOBook}
{Oberguggenberger, M.}
\newblock {\em Multiplication of Distributions and Applications to Partial
  Differential Equations}, volume~{\bf 259} of {\em Pitman Research Notes in
  Mathematics}.
\newblock Longman, Harlow, 1992.
%
%
\bibitem[{Ste}00]{Stdiss}
{Steinbauer, R.}
\newblock {\em Distributional Methods in General Relativity}.
\newblock PhD thesis, University of Vienna, 2000.
%
\end{thebibliography}

\end{document}